\documentclass[letterpaper, 10 pt, conference]{ieeeconf}  

\IEEEoverridecommandlockouts                              

\usepackage{graphics} 
\usepackage{amsmath} 
\usepackage{amssymb}  
\usepackage{tikz}
\usepackage{pgfplots}
\pgfplotsset{width=\columnwidth,compat=1.18} 
\usetikzlibrary{patterns}
\usepgfplotslibrary{fillbetween}

\newtheorem{theorem}{Theorem}
\newtheorem{proposition}{Proposition}
\newtheorem{definition}{Definition}
\newtheorem{observation}{Observation}

\title{\LARGE \bf
Delay effects on the discontinuous stabilization of the nonholonomic integrator and its generalizations*
}

\author{William Clark$^{1}$ and Anthony Bloch$^{3}$
\thanks{*W. Clark was funded by AFOSR grant FA9550-23-1-0400.  A.Bloch was partially supported by NSF grant  DMS-2103026, and AFOSR grants FA
9550-22-1-0215 and FA 9550-23-1-0400.}
\thanks{$^{1}$William Clark is with the Department of Mathematics,
        Ohio University, Athens, OH, 45701, USA
        {\tt\small clarkw3@ohio.edu}}%
\thanks{$^{3}$A. Bloch is with Department of Mathematics, University of Michigan, Ann Arbor, MI 48109, USA. {\tt\small abloch@umich.edu}}
}

\begin{document}

\maketitle
\thispagestyle{empty}
\pagestyle{empty}

\begin{abstract}

The nonholonomic integrator is a famous example in feedback design - although it is small-time locally controllable to the origin, no continuous feedback law exists. Therefore, any stabilizing feedback laws must be either time-varying or discontinuous. A previously studied discontinuous feedback law stabilizes initial conditions lying between two paraboloids and has a sliding mode on the $xy$-plane.

We investigate the effect of introducing delays into this discontinuous feedback law. To a first-order analysis, the lag causes the sliding mode of the $xy$-plane to bifurcate into two switching regions where the resulting dynamics can be interpreted as a hybrid dynamical system with hysteresis. Counterintuitively, the presence of a delay can actually have a positive effect on both the size of the basin of attraction and the convergence rate of the controller.

We also consider the natural generalization of the nonholonomic integrator to higher dimensions. 

\end{abstract}

\section{INTRODUCTION}

There has been much interest in recent decades in finding methods to stabilize systems which fail necessary conditions for stabilization by smooth or 
even continuous feedback. A famous result is: 

\begin{theorem}[Brockett's necessary condition \cite{brockett1983asymptotic}] \label{thm:brockett}
    If the control system $\dot{x}=f(x,u)$ can be locally asymptotically stabilized by means of continuous stationary feedback laws, then the image of $f$ of every neighborhood of $(0,0)\in \mathbb{R}^n\times\mathbb{R}^m$ is a neighborhood of $0\in\mathbb{R}^n$.
\end{theorem}

An important example of a controllable system that fails Theorem \ref{thm:brockett} is the nonholonomic integrator, see \cite{brockett1982control} and \cite{bloch2015nonholonomic}.
Various related results are also discussed in \cite{bloch2015nonholonomic}. Key papers in using time varying control to achieve stabilization include 
\cite{coron1992global} and \cite{pomet1992explicit}. The use of sliding modes is discussed in \cite{BLOCH199691} for the nonholonomic integrator
and more generally in \cite{utkin2004sliding} and \cite{filippov2013differential}. Stabilization of various generalizations of the 
nonholonomic integrator are considered in \cite{bloch2000stabilization}. 

There is also a lot of interest in the effect of delays on the stabilization of various control systems. Interestingly and counterintuitively
it has been shown that delays in linear systems can, under certain conditions, lead to more rapid stabilization, see \cite{moradian2019positive} and 
\cite{qiao2013linear} and related references. A useful recent reference on delays in nonlinear systems is \cite{karafyllis2016recent}.

The contribution of this paper is to present the effect of delays on the sliding mode control for the nonholonomic integrator from \cite{BLOCH199691}. In particular, we show that:
\begin{enumerate}
    \item introducing a delay into a sliding mode causes the mode to bifurcate into two switching surfaces which can be approximated, to first-order, by a hybrid dynamical system,
    \item introducing a delay causes the basin of attraction of the origin to grow, i.e., the feedback control law with delay stabilized more states, and
    \item introducing a delay increases the rate of convergence.
\end{enumerate}

This paper is organized as follows: Section \ref{sec:nh_integrator} introduces the nonholonomic integrator along with the control law and central results from \cite{BLOCH199691}. Section \ref{sec:delayed_controls} introduces the delay in the feedback controls along with a way to compare the infinite-dimensional delay differential equation with the finite-dimensional ordinary differential equation. Section \ref{sec:perturbation} performs a first-order perturbation analysis for small values of delay. This section contains the observation that delays cause the basin of attraction to grow. Section \ref{sec:rate_convergence} compares the rate of convergence against the results in \cite{moradian2019positive}. Section \ref{sec:general} includes a generalization of the nonholonomic integrator with a delay feedback. Finally, conclusions are in Section \ref{sec:conclusion}.


All computations were computed with the DelayDiffEq.jl package in Julia \cite{rackauckas2017differentialequations}.
\section{The Nonholonomic Integrator}\label{sec:nh_integrator}

The nonholonomic ingtegrator equations are given by 
\begin{equation}\label{eq:nh_integrator}
    \dot{x} = u, \quad \dot{y} = v, \quad \dot{z}=xv-yu,
\end{equation}
where $u$ and $v$ are control inputs. Although the system \eqref{eq:nh_integrator} is small-time locally controllable to the origin, there is no solution to
\begin{equation*}
    u = 0, \quad v = 0, \quad xv-yu=\alpha,
\end{equation*}
for $\alpha\ne 0$. Therefore, Theorem \ref{thm:brockett} states that the origin cannot be (locally asymptotically) stabilized by a continuous feedback law.

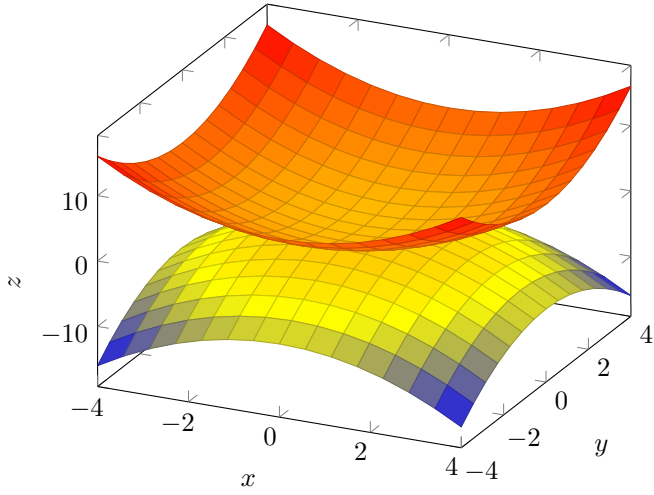
\begin{figure}
    \centering
    \begin{tikzpicture}
        \begin{axis}[xlabel=$x$, ylabel=$y$, zlabel=$z$]
            \addplot3 [surf, samples=15, domain=-4:4]
                {-1/2*(x^2+y^2)};
            \addplot3 [surf, samples=15, domain=-4:4]
                {1/2*(x^2+y^2)};
        \end{axis}
    \end{tikzpicture}
    \caption{Region of initial conditions stabilized by \eqref{eq:no_delay_control} for $\alpha=1$. Initial conditions between the two paraboloids are driven to the origin \eqref{eq:region_no_delay}. No matter how large of a value of $\alpha$ is chosen, the $z$-axis will never lie within $\mathcal{R}_\alpha$.}
    \label{fig:no_delay_region}
\end{figure}

\begin{figure}
	\centering
	\includegraphics[width=\columnwidth]{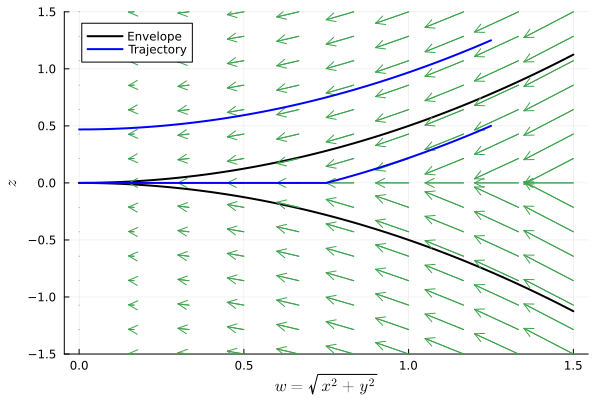}
	\caption{One trajectory within the envelope and one outside. The (green) vector field indicates the control input. The plane $\{z=0\}$ is a sliding mode.}
	\label{fig:no_delay_phase}
\end{figure}

A (discontinuous) control law from \cite{BLOCH199691} has the form
\begin{equation}\label{eq:no_delay_control}
	\begin{split}
		u &= -x + \alpha \cdot y\cdot \mathrm{sign}(z), \\
		v &= -y - \alpha \cdot x\cdot \mathrm{sign}(z),
	\end{split}
\end{equation}
where $\alpha$ is a parameter. Under this control law, the system \eqref{eq:nh_integrator} is Filippov \cite{filippov2013differential} with a sliding mode along the $xy$-plane. Using $V = x^2+y^2$ and $\dot{V} = -2V$, it follows that $x(t)\to 0$ and $y(t)\to 0$ for any initial conditions. The $z$-component approaches zero when
\begin{equation*}
	\frac{1}{2}\alpha (x_0^2 + y_0^2) < |z_0|.
\end{equation*}
Denote the set of all initial conditions obeying the above condition as
\begin{equation}\label{eq:region_no_delay}
    \mathcal{R}_\alpha = \left\{(x,y,z) \in \mathbb{R}^3 : \dfrac{1}{2}\alpha(x^2+y^2) < |z| \right\},
\end{equation}
see Figure \ref{fig:no_delay_region}. If the initial conditions do not lie within $\mathcal{R}_\alpha$, then $z(t)\to z^* \ne 0$, see Figure \ref{fig:no_delay_phase}. 

\section{Delayed Controls}\label{sec:delayed_controls}

We now consider the additional effect of delays.
Consider the lagged version of the control law \eqref{eq:no_delay_control}:
\begin{equation}\label{eq:delayed_controls}
    \begin{split}
        u_\tau &= -x_\tau + \alpha \cdot y_\tau \cdot \mathrm{sign}(z_\tau), \\
        v_\tau &= -y_\tau - \alpha\cdot x_\tau \cdot \mathrm{sign}(z_\tau),
    \end{split}
\end{equation}
where the subscript indicates delay,
\begin{equation*}
    x_\tau(t) = x(t-\tau), \quad y_\tau(t) = y(t-\tau), \quad z_\tau(t) = z(t-\tau),
\end{equation*}
for a fixed constant delay $\tau > 0$.
Applying these controls results in a non-linear delay differential equation:
\begin{equation}\label{eq:delayed_feedback_dynamics}
    \begin{split}
        \dot{x} &=  -x_\tau + \alpha \cdot y_\tau \cdot \mathrm{sign}(z_\tau), \\
        \dot{y} &= -y_\tau - \alpha\cdot x_\tau \cdot \mathrm{sign}(z_\tau),\\
        \dot{z} &= x(-y_\tau - \alpha\cdot x_\tau \cdot \mathrm{sign}(z_\tau)) \\
        &\qquad - y(-x_\tau + \alpha \cdot y_\tau \cdot \mathrm{sign}(z_\tau)).
    \end{split}
\end{equation}
Unlike ordinary differential equations, the initial condition for delay differential equations is an entire history function $\varphi\in C_\tau := C([-\tau,0],\mathbb{R}^3)$ \cite{hale_lunel}, the Banach space of continuous functions. The set of initial conditions \eqref{eq:region_no_delay} becomes
\begin{equation}\label{eq:region_delay}
    \mathcal{R}_{\alpha,\tau} := \left\{ \varphi\in C_\tau : T_t(\varphi)\to 0\right\},
\end{equation}
where $T_t:C_\tau\to C_\tau$ is the solution map.
A direct comparison of $\mathcal{R}_{\alpha,\tau}$ with $\mathcal{R}_\alpha$ is not straightforward as there is not a natural way to embed $\mathbb{R}^3 \hookrightarrow C_\tau$. We will choose an embedding $\iota$ such that
\begin{enumerate}
    \item[(I1)] the initial condition coincides with the last value of the history function, i.e., $\iota(p)(0) = p$ for $p\in\mathbb{R}^3$, and
    \item[(I2)] the solution is differentiable at $t=0$, i.e., if $\gamma(t)$ is a solution with history $\iota(p)$, then
    \begin{equation*}
        \lim_{t\to 0^-} \frac{d}{dt}\iota(p)(t) = \lim_{t\to 0^+} \frac{d}{dt} \gamma(t).
    \end{equation*}
\end{enumerate}
For a given embedding, we will denote 
\begin{equation*}
    \mathcal{R}_{\alpha,\tau}^\iota := \left\{ p\in \mathbb{R}^3 : T_t(\iota(p)) \to 0 \right\}.
\end{equation*}
Notice that for two different embeddings, $\iota \ne \kappa$, we expect that $\mathcal{R}_{\alpha,\tau}^\iota \ne \mathcal{R}_{\alpha,\tau}^\kappa$.

We end this section with a useful observation about the dynamics \eqref{eq:delayed_feedback_dynamics}.
\begin{definition}
    A delay differential equation is equivariant under a group action if
    \begin{equation*}
        f(\gamma .x, \gamma.x_\tau) = \gamma.f(x,x_\tau),
    \end{equation*}
    for all $\gamma\in G$.
\end{definition}
Equivariance of the differential equation carries over to the solution map:
\begin{equation*}
    T_t(\gamma.\varphi) = \gamma.T_t(\varphi).
\end{equation*}

\begin{figure}
    \centering
    \includegraphics[width=\columnwidth]{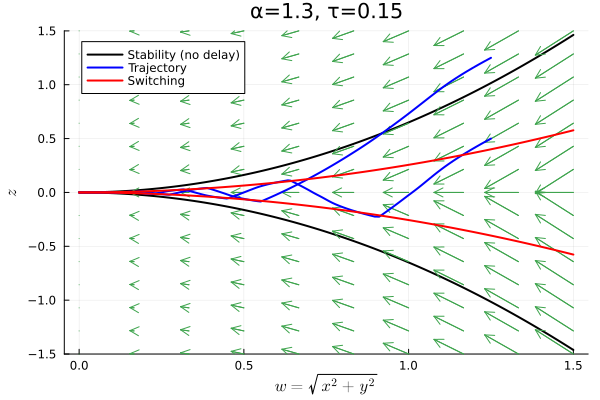}
    \caption{A plot of two trajectories of \eqref{eq:delayed_feedback_dynamics} with one outside of $\mathcal{R}_\alpha$ which indicates that $\mathcal{R}_\alpha \subset \mathcal{R}_{\alpha,\tau}^{\iota^+}$. The two red curves are the boundaries of $P_\pm$ from \eqref{eq:positive_zone} and \eqref{eq:negative_zone}. The region between the two red curves is $P_+\cap P_-$ where hysteresis is observed.}
    \label{fig:delay_phase}
\end{figure}

\begin{proposition}\label{prop:equivariant}
    The dynamics \eqref{eq:delayed_feedback_dynamics} are equivariant under the $\mathrm{SO}_2$ action rotating the $xy$-plane.
\end{proposition}
This proposition allows allows us to set $y(0) = 0$ and analyze $\mathcal{R}_{\alpha,\tau}^\iota$ radially as in Figure \ref{fig:no_delay_phase}, see Figure \ref{fig:delay_phase}.
\subsection{History embedding}
We conclude this section with an explicit choice of history embedding $\iota:\mathbb{R}^3\to C_\tau$ satisfying both (I1) and (I2). For the sake of simplicity, we will take $\iota(p)$ to be linear:
\begin{equation*}
    \iota(p) : \begin{cases}
        x(s) = x_0 + vs, \\
        y(s) = y_0 + ws, \\
        z(s) = z_0 + \nu t.
    \end{cases}
\end{equation*}
To apply the differentiability condition (I2), the sign of $z(-\tau)$ is needed which is independent of the value at $t=0$ (I1). As such, a choice on the sign of $z(-\tau)$ is required which results in two solutions for $v$, $w$, and $\nu$:
\begin{equation*}
    \iota^\pm : \begin{cases}
        \mu_\pm = \dfrac{(\tau\alpha^2-1+\tau)x_0\pm\alpha y_0}{(1-\tau)^2+\alpha^2\tau^2}, \\
        \nu_\pm = \dfrac{\mp\alpha x_0 + (\tau\alpha^2-1+\tau)y_0}{(1-\tau)^2+\alpha^2\tau^2}, \\
        \xi_\pm = x_0 \nu_\pm - y_0 \mu_\pm
    \end{cases}
\end{equation*}
Throughout the remainder of this paper, the embeddings $\iota^\pm$ will be exclusively used. In view of Proposition \ref{prop:equivariant}, $y_0=0$ and $\mathcal{R}_{\alpha,\tau}^{\iota^\pm} \subset \mathbb{R}^2$.
\section{Perturbation analysis for small delays}\label{sec:perturbation}
\begin{figure}
    \centering
    \begin{tikzpicture}
        \draw[->, thick] (-1.5,0) -- (1.5,0);
        \draw[->, thick] (0,-1.5) -- (0,1.5);
        \draw[pattern=north west lines] plot[smooth, samples=100, domain=-1.5:1.5] (\x, \x*\x/1.5) -- (1.5,-1.5) -- (-1.5,-1.5) -- cycle;
        \draw[thick, blue, domain=-1.5:1.5, variable=\x] plot (\x, \x*\x/1.5);
        \node[fill=white] at (-1,-1) {$P_-$};
        \draw[->, thick] (2.5,0) -- (5.5,0);
        \draw[->, thick] (4,-1.5) -- (4,1.5);
        \draw[pattern=north west lines] plot[smooth, samples=100, domain=-1.5:1.5] (\x+4, -\x*\x/1.5) -- (5.5,1.5) -- (2.5,1.5) -- cycle;
        \draw[thick, red, domain=-1.5:1.5, variable=\x] plot (\x+4, -\x*\x/1.5);
        \node[fill=white] at (5,1) {$P_+$};
        \node[blue] at (-0.75, 1) {$\Sigma_-^+$};
        \node[red] at (4.75,-1) {$\Sigma_+^-$};
        \draw[->, thick, blue] (1.25,0.75) -- (2.75, 0.75);
        \draw[->, thick, red] (2.75,-0.75) -- (1.25,-0.75);
        \node[above, blue] at (2,0.75) {$\Delta_-^+$};
        \node[below, red] at (2,-0.75) {$\Delta_+^-$};
    \end{tikzpicture}
    \caption{A hybrid systems schematic of the approximation \eqref{eq:coupled_uv}. In the region $P_+$, the control law \eqref{eq:positive_controls} is implemented while $P_-$ uses \eqref{eq:negative_controls}. The control laws switch when the state hits the sets $\Sigma$.}
    \label{fig:hybrid_schematic}
\end{figure}
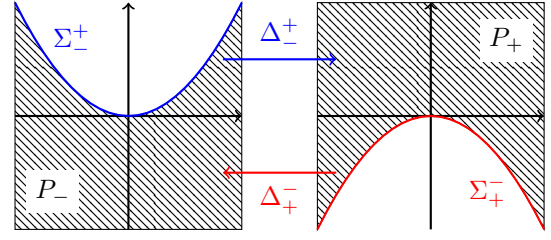

\begin{figure}
    \centering
    \includegraphics[width=\columnwidth]{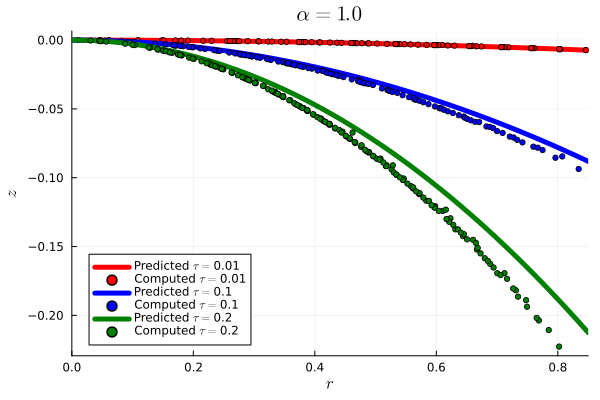}
    \caption{A comparison between the predicted switching curve \eqref{eq:positive_zone} and the location of switching by the actual delay system \eqref{eq:delayed_feedback_dynamics}. The initial conditions are uniformly sampled with $0<x_0,z_0<1$. }
    \label{fig:compare_evelope}
\end{figure}

When $\tau \ll 1$, the DDE \eqref{eq:delayed_feedback_dynamics} can be approximated via 
\begin{equation}\label{eq:first_order_taylor}
    \begin{split}
        x_\tau(t) & \approx x(t) - \tau\dot{x}(t) = x-\tau u, \\
        y_\tau(t) &\approx y(t) -\tau \dot{y}(t) = y-\tau v, \\
        z_\tau(t) &\approx z(t)-\tau\dot{z}(t) = z-\tau(xv-yu).
    \end{split}
\end{equation}
As the approximation \eqref{eq:first_order_taylor} may result in incorrect stability conditions \cite{careful_taylor}, numerical experiments are presented. The controls \eqref{eq:delayed_controls} under the approximation \eqref{eq:first_order_taylor} are
\begin{equation}\label{eq:coupled_uv}
    \begin{split}
        u &= -x + \tau u + \alpha(y-\tau v)\cdot\mathrm{sign}(z-\tau(xv-yu)), \\
        v &= -y + \tau v - \alpha(x-\tau u)\cdot\mathrm{sign}(z-\tau(xv-yu)).
    \end{split}
\end{equation}
Determining the control values requires solving the above coupled set of discontinuous functions \eqref{eq:coupled_uv}. Fortunately, the relationship becomes linear when $z-\tau(xv-yu)$ has constant sign.
\begin{itemize}
    \item Suppose $z>\tau(xv-yu)$. Denote the solutions by $u_+$ and $v_+$:
    \begin{equation}\label{eq:positive_controls}
        \begin{split}
            u_+ &= \frac{(\tau\alpha^2-1+\tau)x + \alpha y}{(1-\tau)^2+\alpha^2\tau^2}, \\
            v_+ &= \frac{-\alpha x + (\tau\alpha^2-1+\tau)y}{(1-\tau)^2+\alpha^2\tau^2}
        \end{split}
    \end{equation}
    This is valid in the region
    \begin{equation}\label{eq:positive_zone}
        \begin{split}
            P_+ &= \left\{ z > \tau(xv_+-yu_+) \right\} \\
            &= \left\{ z > -\dfrac{\alpha\tau}{(1-\tau)^2+\alpha^2\tau^2}(x^2+y^2) \right\}
        \end{split}
    \end{equation}
    \item Suppose $z>\tau(xv-yu)$. Analogously to the other case, we have
    \begin{equation}\label{eq:negative_controls}
        \begin{split}
            u_- &= \frac{(\tau\alpha^2-1+\tau)x - \alpha y}{(1-\tau)^2+\alpha^2\tau^2}, \\
            v_- &= \frac{\alpha x + (\tau\alpha^2-1+\tau)y}{(1-\tau)^2+\alpha^2\tau^2},
        \end{split}
    \end{equation}
    in the region
    \begin{equation}\label{eq:negative_zone}
        P_- = \left\{ z < \dfrac{\alpha\tau}{(1-\tau)^2+\alpha^2\tau^2}(x^2+y^2) \right\}.
    \end{equation}
\end{itemize}
A schematic of these two regions is shown in Figure \ref{fig:hybrid_schematic}.
As $P_+\cap P_-\ne\emptyset$, the approximate system \eqref{eq:first_order_taylor} experiences hysteresis; hence we expect that the $z$-coordinate will decay in oscillatory fashion (Figure \ref{fig:delay_phase}) rather than monotonically (Figure \ref{fig:no_delay_phase}). A comparison between the predicted switching curves and the switching locations from the DDE \eqref{eq:delayed_feedback_dynamics} are shown in Figure \ref{fig:compare_evelope}.

To approximate the region $\mathcal{R}_{\alpha,\tau}$, we notice that the function $V = x^2+y^2$ evolves as
\begin{equation}\label{eq:V_dynamics}
    \dot{V} = 2 \left( \frac{\alpha^2\tau + \tau - 1}{(1-\tau)^2+\alpha^2\tau^2} \right) V.
\end{equation}
The $z$-dynamics in \eqref{eq:delayed_feedback_dynamics} with the approximation \eqref{eq:first_order_taylor} is
\begin{equation}\label{eq:z_dynamics}
    \dot{z} = \frac{-2\alpha}{(1-\tau)^2+\alpha^2\tau^2}V,
\end{equation}
which results in the region
\begin{equation}\label{eq:predicted_boundary}
    \mathcal{R}_{\alpha,\tau}^{\eqref{eq:first_order_taylor}} = \left\{ -\frac{1}{2}\frac{\alpha}{\alpha^2\tau + \tau - 1}(x^2+y^2) < |z| \right\},
\end{equation}
where the superscript indicates that this region comes from the approximation \eqref{eq:first_order_taylor} rather than an embedding.

\begin{figure}
    \centering
    \includegraphics[width=\columnwidth]{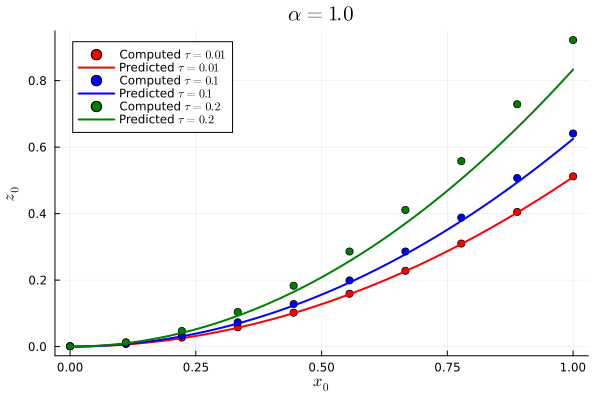}
    \caption{The computed boundary of $\mathcal{R}_{\alpha,\tau}^{\iota^+}$ versus the predicted boundary arising from the approximation \eqref{eq:predicted_boundary}.}
    \label{fig:convergence_boundary}
\end{figure}

\begin{figure}
    \centering
    \includegraphics[width=\columnwidth]{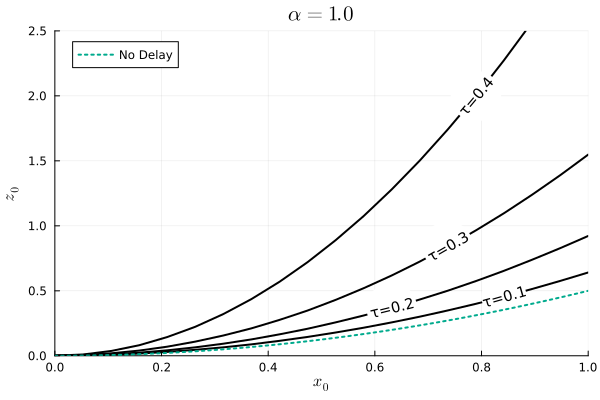}
    \caption{The boundary of $\mathcal{R}_{\alpha,\tau}^{\iota^+}$ for various values of $\tau$ against $\mathcal{R}_{\alpha}$. Notice that the regions grow with $\tau$, $\mathcal{R}_\alpha \subset \mathcal{R}_{\alpha,\tau_1}^{\iota^+} \subset \mathcal{R}_{\alpha,\tau_2}^{\iota^+}$ for (small) $\tau_1 < \tau_2$. When $\tau$ becomes large, the oscillations become large and convergence fails.}
    \label{fig:more_boundaries}
\end{figure}

Notice that the approximation $\mathcal{R}_{\alpha,\tau}^{\eqref{eq:first_order_taylor}}$ was completely agnostic to the choice of embedding with the exception that the solution was assumed to be differentiable, condition (I2). Additionally, when $\tau>0$, $1-\tau - \alpha^2\tau < 1$, which leads to
\begin{equation*}
    \mathcal{R}_\alpha \subset \mathcal{R}_{\alpha,\tau}^{\eqref{eq:first_order_taylor}} \approx \mathcal{R}_{\alpha,\tau}^\iota.
\end{equation*}
\begin{observation}
    Introducing a (small) delay enlarges the basin of attraction of the origin, see Figures \ref{fig:convergence_boundary} and \ref{fig:more_boundaries}.
\end{observation}
\subsection{A hybrid systems interpretation}
The dynamics \eqref{eq:positive_controls} - \eqref{eq:negative_zone} can be viewed as a hybrid dynamical system \cite{hybrid_systems}
\begin{equation}
    \mathcal{H} : \begin{cases}
        \dot{p} = F(p), & p\in C, \\
        p^+ = G(p), & p\in D.
    \end{cases}
\end{equation}
The data is
\begin{equation*}
    \begin{split}
        C &= (P_+\times\{1\}) \cup (P_- \times\{-1\}) \\
        D &= (\Sigma_+^-\times \{1\}) \cup (\Sigma_-^+ \times \{-1\}) \\
        F&= \begin{cases}
            \eqref{eq:nh_integrator}~ \& ~\eqref{eq:positive_controls}, & q=1 \\
            \eqref{eq:nh_integrator}~ \& ~\eqref{eq:negative_controls}, & q=-1
        \end{cases}\\
        G(p, q) &= (p, -q).
    \end{split}
\end{equation*}
See Figure \ref{fig:hybrid_schematic}. This leads to another observation.
\begin{observation}
    A Filippov system with delays can be approximated (to first-order) by a hybrid dynamical system. 
\end{observation}

Moreover, for a given point on a guard $p\in D$, we can compute the time until the guard is reached again. Using \eqref{eq:V_dynamics} and \eqref{eq:z_dynamics}, we get
\begin{equation*}
    t_{\text{switch}} = 2\left(\frac{\kappa}{1-\tau-\alpha^2\tau} \right) \cdot \ln \left( 3 - \frac{2-2\tau}{\kappa} \right),
\end{equation*}
where $\kappa = (1-\tau)^2 + \alpha^2\tau^2$, which is a constant. This leads to our final observation.
\begin{observation}
    Although the trajectory approaches a blocking state, 
    \begin{equation*}
        w \to 0 \in G(D)\cap D,
    \end{equation*}
    the trajectory is not Zeno, i.e. even though Figure \ref{fig:delay_phase} appears to be chattering \cite{fuller} as the trajectory reflects off of the switching curves infinitely many times, it is not.

\end{observation}
\section{Rate of Convergence}\label{sec:rate_convergence}
\begin{figure}
    \centering
    \includegraphics[width=\columnwidth]{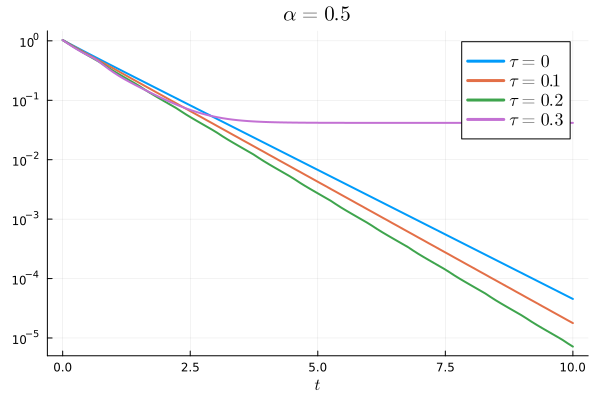}
    \caption{A log plot of $\lVert T_t(\iota^+(p))\rVert$ versus time. As $|\alpha|<1$, the rate of convergence is expected to increase by Theorem \ref{thm:convergence_speedup}. This is true for small enough values of $\tau$. However, when $\tau$ becomes too large, convergence fails altogether. The initial condition is $(x_0,z_0) = (1.0, 0.25)$.}
    \label{fig:trajectory_decay}
\end{figure}

\begin{figure}
    \centering
    \includegraphics[width=\columnwidth]{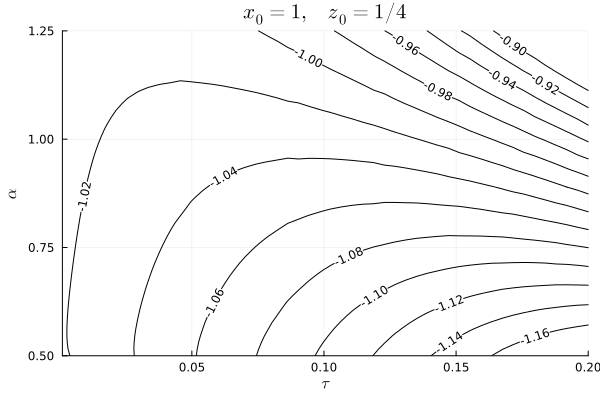}
    \caption{A contour plot of the convergence rate as a function of $\tau$ and $\alpha$. As the initial condition is $(x_0,z_0) = (1.0, 0.25)$, $\alpha > 0.5$ by \eqref{eq:region_no_delay}. Notice that there is still an increase in convergence rate for small delays for some $\alpha > 1$.}
    \label{fig:exponential_contour}
\end{figure}

In Section \ref{sec:perturbation}, it is shown that the introduction of a delay increases the basin of attraction of the origin. A natural follow up question is: How does the introduction of a delay influence the rate of convergence? For the standard first-order example $\dot{x}(t) = -\lambda x(t-\tau)$, the characteristic exponents are given by
\begin{equation*}
    s_k = \frac{1}{\tau} W_k(-\lambda\tau),
\end{equation*}
where $W_k$ is the Lambert-W function. It is straightforward to see that the convergence slows down when $\tau>0$, i.e., $s_0(\tau) > -\lambda$. However, this phenomenon does not hold for all linear systems, \cite{moradian2019positive} and \cite{qiao2013linear}. Consider the linear system of delay differential equations:
\begin{equation*}
    \dot{x}(t) = Ax(t-\tau).
\end{equation*}
\begin{theorem}[Theorem IV.1 in \cite{moradian2019positive}]\label{thm:convergence_speedup}
    Suppose that every eigenvalue of $A$ has negative real part and ordered such that
    \begin{equation*}
        |\mathrm{Re}(\lambda_1)| \leq |\mathrm{Re}(\lambda_2)| \leq \ldots \leq |\mathrm{Re}(\lambda_n)|.
    \end{equation*}
    If for all eigenvalues $\lambda_k$ such that $\mathrm{Re}(\lambda_k) = \mathrm{Re}(\lambda_1)$, we have
    \begin{equation}\label{eq:argument_bound}
        \arg(\lambda_k) \in \left( \frac{3\pi}{4}, \frac{5\pi}{4}\right),
    \end{equation}
    then there exists a delay $\tau\in (0, \bar{\tau})$ for which the delayed rate of convergence is increased $\rho_\tau > \rho_0$, where
    \begin{equation*}
        \begin{split}
            \bar{\tau} &= \min_i \, \frac{1}{|\lambda_i|}\left| \arctan \frac{\mathrm{Re}(\lambda_i)}{\mathrm{Im}(\lambda_i)} \right|, \\
            \rho_0 &= |\mathrm{Re}(\lambda_1)|, \\
            \rho_\tau &= \min_i \, -\frac{1}{\tau} \mathrm{Re}\left (W_0(\lambda_i\tau)\right).
        \end{split}
    \end{equation*}
\end{theorem}
When the $z$-component in \eqref{eq:delayed_feedback_dynamics} has constant sign, the $(x,y)$-dynamics are linear with
\begin{equation*}
    A = \begin{bmatrix}
        -1 & \pm \alpha \\ \mp\alpha & -1
    \end{bmatrix} \implies \lambda_{1,2} = -1 \pm \alpha i.
\end{equation*}
As such, the argument of the eigenvalues are bounded by \eqref{eq:argument_bound} when $-1 < \alpha < 1$. This speedup can be seen in Figure \ref{fig:trajectory_decay} for the specific case of $\alpha=0.5$ and more generally in Figure \ref{fig:exponential_contour}. An increase in convergence rate is still possible when $\alpha>1$.

\section{General Lie algebraic setting}\label{sec:general}
The above analysis for the system \eqref{eq:nh_integrator} with the delayed controls \eqref{eq:delayed_controls} can be carried out for more generalized systems as studied in \cite{bloch2000stabilization} and \cite{Brockett1981}:
\begin{equation*}
    \begin{split}
        \dot{x} &= u, \\
        \dot{Y} &= xu^\top - ux^\top,
    \end{split}
\end{equation*}
where $x,u\in\mathbb{R}^n$ and $Y \in \mathfrak{so}_n$, the special orthogonal Lie algebra consisting of skew-symmetric matrices.
The importance of this system is that it is a canonical form of controllable systems of the form $\dot{x} = B(x)u$ $u\in \mathbb{R}^n$, $x\in \mathbb{R}^{n(n+1)/2}$.
The class in question is the controllable systems of this type where the
first derived algebra of control vector fields spans the
tangent space $T{\mathbb{R}}^{n(n+1)/2}$ at any point. Similarly to \eqref{eq:nh_integrator}, a stabilizing continuous feedback law does not exist.

A different generalization of the Heisenberg system \eqref{eq:nh_integrator}
 is obtained by identifying the vectors $(x,y)^\top$
and $(u,v)^\top$ with the matrices
\begin{equation*}
    X = \dfrac{1}{\sqrt{2}} \begin{bmatrix} x & -y \\ -y & -x \end{bmatrix}, \quad U = \dfrac{1}{\sqrt{2}}\begin{bmatrix} u & -v \\ -v & -u \end{bmatrix},
\end{equation*}
respectively. This yields
\begin{equation*}
    [U,X] = UX-XU = \begin{bmatrix}
        0 & (xv-yu) \\ -(xv-yu) & 0
    \end{bmatrix}.
\end{equation*}
This suggests the following matrix system evolving on $\mathfrak{sl}_n$, the Lie algebra of traceless $n\times n$ matrices:
\begin{equation}\label{eq:sl(n)}
    \dot{X} = U, \quad
    \dot{Y} = [U, X],
\end{equation}
where $X,U \in \mathrm{sym}_{0,n}$ are traceless $n\times n$ real symmetric matrices and $Y\in\mathfrak{so}_n$. Note that
\begin{equation*}
    \mathfrak{sl}_n = \mathrm{sym}_{0,n} \otimes \mathfrak{so}_n,
\end{equation*}
is a direct sum.

More generally, we can consider a Lie algebra with Cartan decomposition
\begin{equation*}
    \mathfrak{g} = \mathfrak{h} \oplus \mathfrak{m},
\end{equation*}
with the flow
\begin{equation}\label{eq:genLie}
    \dot{x} = u, \quad \dot{Y} = [u, x]
\end{equation}
with $x,u\in\mathfrak{m}$ and $Y\in\mathfrak{h}$.

Clearly, the $\mathfrak{sl}_n$ system \eqref{eq:sl(n)} has the form \eqref{eq:genLie}. Additionally, the $\mathfrak{so}_n$ system can also be written in this form. 

Let $\mathfrak{h} = \mathfrak{so}_n$ and $\mathfrak{m} = \mathbb{R}^m$. For $x,u\in\mathfrak{m}$, define $[u,x] := xu^\top - ux^\top \in \mathfrak{h}$. For $Y\in\mathfrak{h}$ and $x\in \mathfrak{m}$, define $[Y,x] = -[x,Y] := Yx$. Then
\begin{equation*}
    \mathfrak{g} := \mathfrak{m} \oplus \mathfrak{h} \cong \mathfrak{so}_{n+1}
\end{equation*}
are isomorphic as Lie algebras through the identification
\begin{equation*}
    \begin{split}
        \mathfrak{h} &\cong \left\{ \begin{bmatrix} 0 & 0 \\ 0 & Y \end{bmatrix} : Y \in \mathfrak{so}_n \right\}, \\
        \mathfrak{m} &\cong \left\{ \begin{bmatrix} 0 & -x^\top \\ x & 0 \end{bmatrix} : x \in \mathbb{R}^n \right\}.
    \end{split}
\end{equation*}
A (temporally) discontinuous controller for these systems is given in \cite{bloch2000stabilization}.

\subsection{Three dimensions}\label{sec:3_example}
For the case when $n=3$, we have the identification $(\mathfrak{so}_n, [\cdot,\cdot]) \cong (\mathbb{R}^3, \times)$, where $\times$ is the cross-product. The corresponding control system is
\begin{equation*}
    \dot{x} = u, \quad \dot{y} = x\times u,
\end{equation*}
where $x, u, y\in\mathbb{R}^3$. An analogous feedback law to \eqref{eq:no_delay_control} is
\begin{equation}\label{eq:so3_control}
    u = -x + \alpha \cdot \tilde{y} \times x,
\end{equation}
where
\begin{equation*}
    \tilde{y} = \begin{bmatrix}
        \mathrm{sign}(y_1) & \mathrm{sign}(y_2) & \mathrm{sign}(y_3)
    \end{bmatrix}^\top.
\end{equation*}
Using the function $V = \lVert x \rVert^2$, we have
\begin{equation*}
    \begin{split}
        \dot{V} &= 2\langle x, \dot{x}\rangle = -2V, \\
        \dot{y} &= \alpha \cdot x\times (\tilde{y}\times x) = \alpha\cdot \left[ V\cdot \tilde{y} - \langle x, \tilde{y}\rangle \cdot x \right].
    \end{split}
\end{equation*}
As a result, $x\to 0$. However, the analysis for $y$ is more complicated as another term would appear in \eqref{eq:z_dynamics}.

\begin{figure}
    \centering
    \includegraphics[width=\columnwidth]{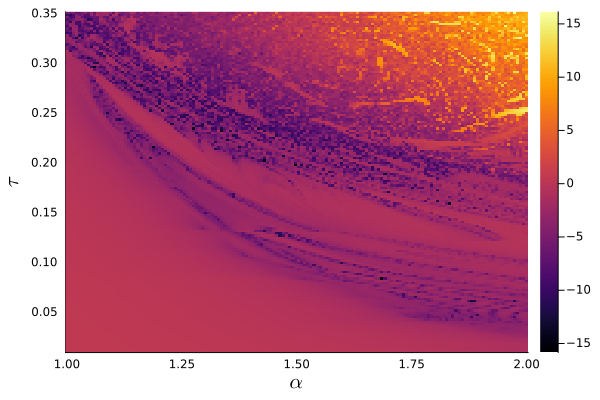}
    \caption{A heatmap of $\log \lVert y(20)\rVert$ against $(\alpha,\tau)$ for the example in Section \ref{sec:3_example}. While the state $y$ fails to be driven to the origin when $\tau=0$, there are certain values of $\tau>0$ where $y\to 0$.}
    \label{fig:so_3_heatmap}
\end{figure}

In this case, the control law \eqref{eq:so3_control} has sliding modes on all coordinate planes. Introducing a (small) delay will cause the planes to bifurcate into families of paraboloids. A convergence plot for various values of $\alpha$ and $\tau$ corresponding to the initial condition
\begin{equation*}
    x = \begin{bmatrix} 0.2 & 1.1 & 1.1 \end{bmatrix}^\top, \quad y = \begin{bmatrix} 3.0 & -0.2 & 0.1 \end{bmatrix}^\top,
\end{equation*}
are shown in Figure \ref{fig:so_3_heatmap}. This problem and initial conditions match Section 6 in \cite{bloch2000stabilization}.

\section{CONCLUSIONS}\label{sec:conclusion}

In this paper we considered non-smooth stabilization of the nonholonomic integrator in the presence of delays. 
This system is a classic example where continuous stabilization to the origin is impossible despite being controllable. Through a combination of perturbation analysis, spectral analysis, and numerical analysis, we showed that introducing a delay can stabilize more states and increase the rate of convergence. Finally, we examine a generalization to the rotation Lie algebra $\mathfrak{so}_3$.

In the future, we intend to analyze the extension to general Lie algebras as presented in Section \ref{sec:general}. In particular, we intend to provide perturbative analysis to the control system \eqref{eq:genLie} with delayed feedback.

\addtolength{\textheight}{-12cm}   









\bibliographystyle{IEEEtran}
\bibliography{references}

\end{document}